# First applications of generalized Li's criterion to study the Riemann zeta-function zeroes location

**S. K. Sekatskii** (LPMV, Ecole Polytechnique Fédérale de Lausanne, Switzerland)

We present the first applications of the recently established by us (arXiv:1304.7895; *Ukrainian Math. J.* – 2014.- **66**. - P. 371 – 383) generalized Li's criterion equivalent to the Riemann Hypothesis. This criterion is the statement that the Riemann hypothesis is equivalent to the non-negativity of the derivatives $\frac{1}{(m-1)!}\frac{d^m}{dz^m}((z+b)^{m-1}\ln(\xi(z)))|_{z=1+b}$ of the Riemann xi-function for all real $b>-1/2$ and all $m$=1, 2, 3… We show that for any positive integer $n$ there is such value of $b_n$ (depending on $n$) that for all $m \le n$ and $b>b_n$, inequality $\frac{1}{(m-1)!}\frac{d^m}{dz^m}((z+b)^{m-1}\ln(\xi(z)))|_{z=1+b} \ge 0$ does hold true. Assuming RH, we also have found an asymptotic of the generalized Li's sums over non-trivial Riemann zeroes for large $n$, and discuss what asymptotic of $\frac{1}{(n-1)!}\frac{d^n}{dz^n}((z+b)^{n-1}\ln((z-1)\varsigma(z)))|_{z=1+b}$ is required for the Riemann hypothesis holds true.

## Introduction

Recently, in Ref. [1] we have established the generalized Bombieri – Lagarias' theorem (see [2] for their original theorem) and the generalized Li's criterion of the truth of the Riemann hypothesis concerning the location of non-trivial zeroes of the Riemann zeta-function (see [3] for the original Li's criterion, and see e.g. [4] for standard definitions and discussion of the general properties of the Riemann zeta-function):

**Theorem 1. (Generalized Bombieri – Lagarias' theorem).** *Let a and* $\sigma$ *are arbitrary real numbers,* $a < \sigma$, *and R be a multiset of complex numbers* $\rho$ *such that*

*(i)* $2\sigma - a \notin R$

*(ii)* $\sum_{\rho} (1 + |\operatorname{Re}\rho|)/(1 + |\rho + a - 2\sigma|^2) < +\infty$

*Then the following conditions are equivalent*

*(a)* $\operatorname{Re}\rho \le \sigma$ *for every* $\rho$*;*

*(b)* $\sum_{\rho} \operatorname{Re}(1 - \left(\frac{\rho - a}{\rho - 2\sigma + a}\right)^n) \ge 0$ *for n=1, 2, 3...*

*(c) For every fixed* $\varepsilon > 0$ *there is a positive constant* $c(\varepsilon)$ *such that* $\sum_{\rho} \operatorname{Re}(1 - \left(\frac{\rho - a}{\rho - 2\sigma + a}\right)^n) \ge -c(\varepsilon)e^{\varepsilon n}$*, n=1, 2, 3...*

*If at the same conditions* $a > \sigma$ *is taken, the point (a) is to be changed to*

*(a')* $\operatorname{Re}\rho \ge \sigma$ *for every* $\rho$*,*

*points (b), (c) remain unchanged.*

*If, additionally to the aforementioned conditions, also the following takes place:*

*(iii) If* $\rho \in R$*, than* $\bar{\rho} \in R$ *with the same multiplicity as* $\rho$

*one can omit the operation of taking the real part in (b), (c), the expressions at question are real. (Here, as usual, $\bar{\rho}$ means a complex conjugate of $\rho$).*

**Theorem 2. (Generalized Li's criterion).** *Let a is an arbitrary real number, $a \neq \sigma$, and R be a multiset of complex numbers $\rho$ such that*

*(i)* $2\sigma - a \notin R$, $a \notin R$

*(ii)* $\sum_{\rho}(1+|\mathrm{Re}\,\rho|)/(1+|\rho+a-2\sigma|^2) < +\infty$, $\sum_{\rho}(1+|\mathrm{Re}\,\rho|)/(1+|\rho-a|^2) < +\infty$

*(iii)* *If* $\rho \in R$, *than* $2\sigma - \rho \in R$

*Then the following conditions are equivalent*

*(a)* $\mathrm{Re}\,\rho = \sigma$ *for every* $\rho$;

*(b)* $\sum_{\rho}\mathrm{Re}(1-\left(\frac{\rho-a}{\rho+a-2\sigma}\right)^n) \geq 0$ *for any a and n=1, 2, 3...*

*(c) For every fixed* $\varepsilon > 0$ *and any a there is a positive constant* $c(\varepsilon, a)$ *such that* $\sum_{\rho}\mathrm{Re}(1-\left(\frac{\rho-a}{\rho+a-2\sigma}\right)^n) \geq -c(\varepsilon, a)e^{\varepsilon n}$, *for n=1, 2, 3...*

*If, additionally to the aforementioned conditions, also the following takes place:*

*(iv)* *If* $\rho \in R$, *than complex conjugate* $\bar{\rho} \in R$ *with the same multiplicity as* $\rho$

*one can omit the operation of taking the real part in (b), (c), the expressions at question are real.*

Then, applying the generalized Littlewood theorem about contour integrals of logarithm of an analytical function (see below and [1, 5-7]), we have established the equality

$$\frac{1}{n(2a-1)}\sum_{\rho}(1-\left(\frac{\rho-a}{\rho+a-1}\right)^n) = \frac{1}{(n-1)!}\frac{d^n}{dz^n}((z-a)^{n-1}\ln(\xi(z)))|_{z=1-a}$$

and have proven also the following [1]:

**Theorem 3.** *Riemann hypothesis is equivalent to the non-negativity of all derivatives* $\frac{1}{(n-1)!}\frac{d^n}{dz^n}((z-a)^{n-1}\ln(\xi(z)))|_{z=1-a}$ *for all positive integers n and any real a<1/2; correspondingly, it is equivalent also to the non-positivity of all derivatives* $\frac{1}{(n-1)!}\frac{d^n}{dz^n}((z-a)^{n-1}\ln(\xi(z)))|_{z=1-a}$ *for all positive integers n and any real a>1/2.*

Thus to judge the truth of the Riemann hypothesis, certain derivatives of the Riemann xi-function can be estimated at an arbitrary point of the real axis except the point $z$=1/2, not only at the point $z$=1 (or 0) as this was initially formulated by Li [3]. In particular, this can be done far to the right from the point $z$=1, where Riemann zeta-function and its logarithm are defined by absolutely convergent series [4]

$$\varsigma(z) = \sum_{n=1}^{\infty}\frac{1}{n^z} \tag{1}$$

and

$$\ln\varsigma(z) = -\sum_p \ln\left(1-\frac{1}{p^z}\right) = \sum_p\sum_{n=1}^{\infty}\frac{1}{np^{nz}} = \sum_{n=2}^{\infty}\frac{\Lambda(n)}{\ln n\cdot n^z} \qquad (2).$$

(in (2) we have a sum over primes or use the van Mandgoldt function). This circumstance holds promise to elucidate certain properties of the Riemann function zeroes, and in the current paper we present the first applications of such an approach.

For completeness, we finish the Introduction section presenting the generalized Littlewood theorem.

**Theorem 4 (Generalized Littlewood theorem).** *Let C denotes the rectangle bounded by the lines* $x = X_1,\ x = X_2,\ y = Y_1,\ y = Y_2$ *where* $X_1 < X_2,\ Y_1 < Y_2$ *and let f(z) be analytic and non-zero on C and meromorphic* *inside it, let also g(z) is analytic on C and meromorphic inside it. Let F(z)=ln(f(z)), the logarithm being defined as follows: we start with a particular determination on* $x = X_2$*, and obtain the value at other points by continuous variation along y=const from* $\ln(X_2 + iy)$*. If, however, this path would cross a zero or pole of f(z), we take F(z) to be* $F(z \pm i0)$ *according as we approach the path from above or below. Let also the poles and zeroes of the functions f(z), g(z) do not coincide.*

*Then* $$\int_C F(z)g(z)dz = 2\pi i\left(\sum_{\rho_g} res(g(\rho_g)\cdot F(\rho_g)) - \sum_{\rho_f^0} \int_{X_1+iY_\rho^0}^{X_\rho^0+iY_\rho^0} g(z)dz + \sum_{\rho_f^{pol}} \int_{X_1+iY_\rho^{pol}}^{X_\rho^{pol}+iY_\rho^{pol}} g(z)dz\right)$$

*where the sum is over all* $\rho_g$ *which are poles of the function g(z) lying inside C, all* $\rho_f^0 = X_\rho^0 + iY_\rho^0$ *which are zeroes of the function f(z) counted taking into account their multiplicities (that is the corresponding term is multiplied by m for a zero of the order m) and which lye inside C, and all* $\rho_f^{pol} = X_\rho^{pol} + iY_\rho^{pol}$ *which are poles of the function f(z) counted taking into account their multiplicities and which lye inside C. For this is true all relevant integrals in the right hand side of the equality should exist.*

**2. Relation between the location of the non-trivial Riemann zeta-function zeroes and certain derivatives involving logarithm of this function**

Let us use theorem 3 in the form that for *b*>-1/2 all derivatives $\frac{1}{(n-1)!}\frac{d^n}{dz^n}((z+b)^{n-1}\ln(\xi(z)))|_{z=1+b}$ should be non-negative for the Rieman hypothesis is true, and the known relation between xi- and zeta- Riemann functions:

$$\xi(z)=\frac{1}{2}z(z-1)\pi^{-z/2}\Gamma(z/2)\varsigma(z) \qquad (3).$$

Clearly, we have

$$\frac{d^n}{dz^n}((z+b)^{n-1}\ln(\xi(z)))|_{z=1+b}=\frac{d^n}{dz^n}((z+b)^{n-1}\ln 1/2)|_{z=1+b}+\frac{d^n}{dz^n}((z+b)^{n-1}\ln(z-1))|_{z=1+b}$$

$$+\frac{d^n}{dz^n}((z+b)^{n-1}\ln\pi^{-z/2})|_{z=1+b}+\frac{d^n}{dz^n}((z+b)^{n-1}\ln(z\Gamma(z/2))|_{z=1+b} \qquad (4)$$

$$+\frac{d^n}{dz^n}((z+b)^{n-1}\ln(\varsigma(z)))|_{z=1+b},$$

and now we analyze all terms here one by one.

Trivially $\frac{d^n}{dz^n}((z+b)^{n-1}\ln(1/2))|_{z=1+b}=0$ and $\frac{d^n}{dz^n}((z+b)^{n-1}\ln\pi^{-z/2})|_{z=1+b}=-\frac{\ln\pi}{2}n!$, while to calculate all other terms we apply the generalized Littlewood theorem in a manner similar to that used in our Ref. [1]. Namely, we introduce the function $g(z)=\frac{n(2b+1)(z+b)^{n-1}}{(z-b-1)^{n+1}}-\frac{n(2b+1)}{(z-b-1)^2}$ and rectangular contour $C$ with vertices at $\pm X\pm iX$ with real $X\to+\infty$, if some pole of the gamma – function occurs on the contour just shift it a bit to avoid this, and consider a contour integral $\int_C g(z)\ln(z\Gamma(z/2))dz$. Known asymptotic of the logarithm of the gamma-function for large $|z|$, $\ln\Gamma(z)\cong O(z\ln z)$ guaranties the "disappearance" of the contour integral value (it tends to zero when $X\to\infty$ due to the asymptotic $g(z)\cong O(1/z^3)$) thus we get, after division by $2\pi i$ and elementary transformations

$$\frac{1}{(n-1)!}\frac{d^n}{dz^n}((z+b)^{n-1}\ln(z\Gamma(z/2)))|_{z=b+1}=$$
$$-\frac{1}{2b+1}\sum_{k=1}^{\infty}(1-\left(\frac{2k-b}{2k+b+1}\right)^n-\frac{n(2b+1)}{2k+b+1})+\frac{n}{2}\psi\left(\frac{b+1}{2}\right)+\frac{n}{b+1} \qquad (5).$$

Here $\psi(z)$ is a digamma function. In the sum occurring in this expression one easily recognizes the sum $\sum_{\rho_f^{pol}} \int_{X_1+iY_\rho^{pol}}^{X_\rho^{pol}+iY_\rho^{pol}} g(z)dz$ taken over simple poles of the gamma function (pole at $z$=0 is cancelled by the factor $z$): clearly, $\frac{n(2b+1)(z+b)^{n-1}}{(z-b-1)^{n+1}} = \frac{d}{dz}(1-\left(\frac{z+b}{z-b-1}\right)^n)$ which explains why function *g(z)* is used here; the term $-\frac{n(2b+1)}{(z-b-1)^2}$ is added just to ensure the asymptotic $g(z) \cong O(1/z^3)$ necessary to bring the contour integral value to zero (summands $\frac{n}{2}\psi\left(\frac{b+1}{2}\right)+\frac{n}{b+1}$ in eq. (5) come exactly from this term).

Similarly, using the function $\tilde{g}(z) = \frac{n(2b+1)(z+b)^{n-1}}{(z-b-1)^{n+1}}$ and the same contour, generalized Littlewood theorem gives

$$\frac{1}{(n-1)!}\frac{d^n}{dz^n}((z+b)^{n-1}\ln(z-1))|_{z=b+1} = \frac{1}{2b+1}(1+(-1)^{n+1}(1+\frac{1}{b})^n) \qquad (6).$$

Collecting everything together, we have

$$\frac{1}{(n-1)!}\frac{d^n}{dz^n}((z+b)^{n-1}\ln(\xi(z)))|_{z=1+b} =$$
$$-\frac{1}{2b+1}\sum_{k=1}^{\infty}(1-\left(\frac{2k-b}{2k+b+1}\right)^n - \frac{n(2b+1)}{2k+b+1}) + \frac{n}{2}\psi\left(\frac{b+1}{2}\right) + \frac{n}{b+1} - \frac{n\ln\pi}{2}$$
$$+\frac{1}{2b+1}(1+(-1)^{n+1}(1+\frac{1}{b})^n) + \frac{1}{(n-1)!}\frac{d^n}{dz^n}((z+b)^{n-1}\ln(\varsigma(z)))|_{z=1+b} \qquad (7).$$

The value of $\frac{1}{(n-1)!}\frac{d^n}{dz^n}((z+b)^{n-1}\ln(\varsigma(z)))|_{z=1+b}$ can be expressed via certain integral using the same hint as above. For this we consider a contour integral $\int_C \tilde{g}(z)\ln(\varsigma(z))dz$ with the rectangular contour having the vertices

$c \pm iX$, $c + X \pm iX$ with real $X \to +\infty$ and real $c$, $1 < c < b+1$. Cauchy residue theorem gives

$$\frac{1}{(n-1)!}\frac{d^n}{dz^n}((z+b)^{n-1}\ln(\varsigma(z)))|_{z=1+b} = -\frac{n}{2\pi}\int_{-\infty}^{\infty}\ln(\varsigma(c+it))\frac{(c+b+it)^{n-1}}{(c-b-1+it)^{n+1}}dt \quad (8).$$

Similarly,

$$\frac{1}{(n-1)!}\frac{d^n}{dz^n}((z+b)^{n-1}\ln((z-1)\cdot\varsigma(z)))|_{z=1+b} = -\frac{n}{2\pi}\int_{-\infty}^{\infty}\ln((c-1+it)\cdot\varsigma(c+it))\frac{(c+b+it)^{n-1}}{(c-b-1+it)^{n+1}}dt$$
(8a).

For $\frac{1}{(n-1)!}\frac{d^n}{dz^n}((z+b)^{n-1}\ln(\xi(z)))|_{z=1+b}$ we obtain a bit more complex expression due to the necessity to use the function $g_1(z) = \frac{n(z+b)^{n-1}}{(z-b-1)^{n+1}} - \frac{n}{(z-b-1)^2}$ now:

$$\frac{1}{(n-1)!}\frac{d^n}{dz^n}((z+b)^{n-1}\ln(\xi(z)))|_{z=1+b} - \frac{n\xi'}{\xi}(b+1) = -\frac{n}{2\pi}\int_{-\infty}^{\infty}\ln(\xi(c+it))\frac{(c+b+it)^{n-1}}{(c-b-1+it)^{n+1}}dt$$

One more representation of the derivatives at question deserves a special Lemma.

***Lemma 1.*** For Rez>1,

$$\frac{1}{(n-1)!}\frac{d^n}{dz^n}((z+b)^n\ln(\varsigma(z)) = -\sum_{m=1}^{\infty}\frac{\Lambda(m)}{m^z}L_{n-1}^1((z+b)\cdot\ln m) \qquad (9),$$

where $L_{n-1}^1(x) = -\frac{d}{dx}L_n(x) = \sum_{j=1}^{n}C_n^j\frac{(-1)^j x^{j-1}}{(j-1)!}$ is generalized Laguerre polynomial.

*Proof.* Using Leibnitz rule, we write

$\frac{1}{(n-1)!}\frac{d^n}{dz^n}((z+b)^{n-1}\ln(\varsigma(z))) = \frac{1}{(n-1)!}\sum_{j=1}^{n}C_n^j\frac{d^j(\ln(\varsigma(z))}{dz^j}\cdot\frac{d^{n-j}}{dz^{n-j}}(z+b)^{n-1}$. (The term corresponding to *j=0* is clearly zero). Further,

$$\frac{1}{(n-1)!}\sum_{j=1}^{n}C_n^j\frac{d^j(\ln(\varsigma(z))}{dz^j}\cdot\frac{d^{n-j}}{dz^{n-j}}(z+b)^{n-1} =$$
$$\frac{1}{(n-1)!}\sum_{j=1}^{n}C_n^j\frac{d^{j-1}(\varsigma'(z)/\varsigma(z))}{dz^{j-1}}(n-1)(n-2)...j(z+b)^{j-1} =$$

$\sum_{j=1}^{n} C_n^j \frac{1}{(j-1)!} \frac{d^{j-1}(\varsigma'(z)/\varsigma(z))}{dz^{j-1}}(z+b)^{j-1}$. Using $\frac{\varsigma'(z)}{\varsigma(z)} = -\sum_{m=1}^{\infty} \frac{\Lambda(m)}{m^z}$ valid for Rez>1 [4], we have

$\sum_{j=1}^{n} C_n^j \frac{1}{(j-1)!} \frac{d^{j-1}(\varsigma'(z)/\varsigma(z))}{dz^{j-1}}(z+b)^{j-1} = -\sum_{j=1}^{n} C_n^j \frac{(z+b)^{j-1}}{(j-1)!} \sum_{m=1}^{\infty} \frac{(-1)^{j-1}\Lambda(m)\ln^{j-1} m}{m^z}$ and now we can change the summation order (for the case at hand this operation is certainly legitimate):

$$-\sum_{j=1}^{n} C_n^j \frac{(b+z)^{j-1}}{(j-1)!} \sum_{m=1}^{\infty} \frac{(-1)^{j-1}\Lambda(m)\ln^{j-1} m}{m^z} = -\sum_{m=1}^{\infty} \frac{\Lambda(m)}{m^z} \sum_{j=1}^{n} C_n^j \frac{(b+z)^{j-1}(-1)^{j-1}\ln^{j-1} m}{(j-1)!}.$$

Reminding the definition of Laguerre polynomials, $\sum_{j=0}^{n} C_n^j \frac{(-1)^j x^j}{j!} = L_n(x)$ and $\frac{d}{dx} L_n(x) = \sum_{j=1}^{n} C_n^j \frac{(-1)^j x^{j-1}}{(j-1)!} = -L_{n-1}^1(x)$, we finish the proof of (9).

Evident consequence of the above Lemma is, for $b$>0,

$$\frac{1}{(n-1)!} \frac{d^n}{dz^n}((z+b)^n \ln(\varsigma(z))|_{z=b+1} = -\sum_{m=1}^{\infty} \frac{\Lambda(m)}{m^{b+1}} L_{n-1}^1((2b+1)\ln m) \qquad (10).$$

In relation with these formulae, it is worthwhile to note that Laguerre polynomials $L_{n-1}^1(x)$ have already appeared in Riemann function researches [8 - 10].

**Remark 1.** Application of the same approach to the derivative $\frac{1}{(n-1)!} \frac{d^n}{dz^n}((z+b)^{n-1} \ln(z\Gamma(z/2)))|_{z=b+1}$ without difficulties restores equation (5).

3. **For any *n*, the derivatives at question are indeed positive in the limit of large *b***

Our next aim is to prove the following

**Theorem 5**. *For any arbitrary large integer m there exists positive real c depending on m such that the inequality* $\frac{d^n}{dz^n}((z+b)^{n-1}\ln(\xi(z)))|_{z=1+b} > 0$ *does hold true for all* $n \le m$ *and all* $b \ge c$.

For this we first prove two technical Lemmas.

***Lemma 2***. For any fixed $b \ge 0$ we have in the limit of large *n*:

$$\frac{2b+1}{2}(n\ln n + n\gamma - n + n\ln(2-\frac{5}{b+3}) + O(1)) \le \sum_{k=1}^{\infty}(\left(\frac{2k-b}{2k+b+1}\right)^n - 1 + \frac{n(2b+1)}{2k+b+1}) \le$$
$$\frac{2b+1}{2}(n\ln n + n\gamma - n + n\ln(2-\frac{1}{b+1}) + O(1)) \qquad (11),$$

where $\gamma$ is Euler – Mascheroni constant.

*Proof*. Clearly, $\frac{2k-b}{2k+b+1} = 1 - \frac{2b+1}{2k+b+1}$ and we need to estimate a sum $S = (\sum_{k=1}^{\infty}\left(1-\frac{2b+1}{2k+b+1}\right)^n - 1 + \frac{n(2b+1)}{2k+b+1})$. The function *f(k)* to be summed is a monotonically decreasing function of *k*, so we have $\int_1^{\infty}((1-\frac{2b+1}{2x+b+1})^n - 1 + \frac{n(2b+1)}{2x+b+1})dx \le S \le \int_1^{\infty}((1-\frac{2b+1}{2x+b-1})^n - 1 + \frac{n(2b+1)}{2x+b-1})dx$. We make a variable change $y = \frac{2b+1}{2x+b-1}$ in the right integral and $y = \frac{2b+1}{2x+b+1}$ in the left one obtaining

$$\frac{2b+1}{2}\int_0^{\frac{2b+1}{b+3}}((1-y)^n - 1 + ny)\frac{dy}{y^2} \le S \le \frac{2b+1}{2}\int_0^{\frac{2b+1}{b+1}}((1-y)^n - 1 + ny)\frac{dy}{y^2}.$$

Now, in the integral $\frac{2b+1}{2}\int_0^{\frac{2b+1}{b+1}}((1-y)^n - 1 + ny)\frac{dy}{y^2}$ we divide the integration range into segments $[0,\ 1]$ and $\left[1,\ \frac{2b+1}{b+1}\right]$. For the second segment, the estimations $|\int_1^{\frac{2b+1}{b+1}} \frac{1}{y^2}(1-y)^n dy| < \int_1^2 (y-1)^n dy = \frac{1}{n+1}$, $\int_1^{\frac{2b+1}{b+1}} \frac{1}{y^2} dy < \frac{1}{2}$, $\int_1^{\frac{2b+1}{b+1}} \frac{n}{y} dy = n\ln(2-\frac{1}{b+1})$ suffice. For the first one, we have $\int_0^1 \frac{(1-y)^n - 1 + ny}{y^2} dy = n(\psi(n)+\gamma-1)+1$, as can be seen integrating by parts with *d(1/y)* and applying entry 3.231.5 of GR book [11]: for $\mathrm{Re}\,\nu > 0,\ \mathrm{Re}\,\mu > 0$, $\int_0^1 \frac{x^{\mu-1} - x^{\nu-1}}{1-x} dx = \psi(\nu) - \psi(\mu)$. (This is worthwhile to notice an interesting discussion of similar sums and integrals in Appendix A of [8]). Thus $S \le \frac{2b+1}{2}(n\ln n + n\gamma - n + n\ln 2 + O(1))$; properties $\psi(n) = \ln n - \frac{1}{2n} + O(1/n^2)$ and $\psi(1) = -\gamma$ are used here. The integral $\frac{2b+1}{2}\int_1^{\frac{2b+1}{b+3}} ((1-y)^n - 1 + ny)\frac{dy}{y^2}$ is estimated in a similar fashion, and with the *O(n)* precision, in the segment $\left[1,\ \frac{2b+1}{b+3}\right]$ only the $\int_1^{\frac{2b+1}{b+3}} \frac{n}{y} dy = n\ln(2-\frac{5}{b+3})$ term is important. We have thus $S \ge \frac{2b+1}{2}(n\ln n + n\gamma - n + n\ln(2-\frac{5}{b+3}) + O(1))$ which proves what we want.

***Lemma 3***. For any real *b*>1/2

$$I = \int_{-\infty}^{\infty} \frac{((2b)^2 + t^2)^{(n-1)/2}}{(1+t^2)^{(n+1)/2}} dt \le \frac{2b}{4b^2-1} \cdot (2b)^{n+1} \qquad (12).$$

*Proof.* Let *n-1=2k* even hence we have $I=\int_{-\infty}^{\infty}\frac{((2b)^2+t^2)^k}{(1+t^2)^{k+1}}dt$. First we use the example N 3.251.4 from Gradshtein and Ryzhik book [11]: $\int_{-\infty}^{\infty}\frac{x^{2m}}{(1+x^2)^n}dx=\frac{(2m-1)!!(2n-2m-3)!!\pi}{(2n-2)!!}$. Thus

$I=\sum_{l=0}^{k}C_k^l\{(2b)^2\}^l\int_{-\infty}^{\infty}\frac{t^{2(k-l)}}{(1+t^2)^{k+1}}=\sum_{l=0}^{k}C_k^l\frac{(2k-2l-1)!!(2l-1)!!\pi}{(2k)!!}\{(2b)^2\}^l$. Now we rewrite as $(2n-1)!!=\frac{2^n}{\sqrt{\pi}}\Gamma(n+1/2)$, $(2n)!!=2^n\Gamma(n+1)$ [11] and obtain $I=\sum_{l=0}^{k}C_k^l\frac{\Gamma(k-l+1/2)\Gamma(l+1/2)}{\Gamma(k+1)}\{(2b)^2\}^l$. Reminding that $C_k^l=\frac{\Gamma(k+1)}{\Gamma(l+1)\Gamma(k-l+1)}$, we get $I=\sum_{l=0}^{k}\frac{\Gamma(k-l+1/2)\Gamma(l+1/2)}{\Gamma(l+1)\Gamma(k-l+1)}\{(2b)^2\}^l$. For all *k*, *l*; $l\le k$, we have $\frac{\Gamma(k-l+1/2)}{\Gamma(k-l+1)}<1$ and $\frac{\Gamma(l+1/2)}{\Gamma(l+1)}<1$ so that the coefficients of this expansion are all *<1*. Thus we have $I\le\sum_{l=0}^{k}\{(2b)^2\}^l$ which is a sum of a geometrical progression terms, so that $I\le\frac{(2b)^{2k+2}-1}{4b^2-1}<\approx{C}(2b)^{2k+2}=\approx{C}(2b)^{n+1}$ where $\approx{C}=\frac{1}{4b^2-1}$.

If *n-1=2k+1* is odd, write $\frac{((2b)^2+t^2)^k}{(1+t^2)^k}\cdot\frac{((2b)^2+t^2)^{1/2}}{(1+t^2)^{1/2}}\le 2b\frac{((2b)^2+t^2)^k}{(1+t^2)^k}$ and repeat the same calculations again. Thus we always have $I=\int_{-\infty}^{\infty}\frac{((2b)^2+t^2)^{(n-1)/2}}{(1+t^2)^{(n+1)/2}}dt\le\tilde{C}(2b)^{n+1}$ where $\tilde{C}=\frac{2b}{4b^2-1}$.

***Remark 2.*** With the variable changes $x=(2b)^2+t^2$ and then $x=a/y$ we have for the same integral $I=a^{-3/2}\int_0^1 y^{-1/2}(1-y)^{-1/2}(1-yu/a)^{-(n+1)/2}dy$, where $a=4b^2$, $u=4b^2-1$. This integral is nothing else than a particular case of well-known integral representation of ${}_2F_1\equiv F(\alpha,\ \beta;\ \gamma,\ z)$ hypergeometric function, see e.g. [12, 13]: $\int_0^1 y^{\lambda-1}(1-y)^{\mu-1}(1-\beta y)^{-\nu}\,dy=B(\lambda,\mu)\times F(\lambda,\nu;\lambda+\mu;\beta)$ (adopting for our case, we have $I=\frac{\pi}{8b^3}F(\frac{1}{2},\ \frac{n+1}{2};\ 1;\ 1-\frac{1}{4b^2})$), hence much studied properties of this function apparently might be used for its estimation. Amusingly, the situation is not so simple. For example, on p. 77 of [12], asymptotic of $F(\alpha,\ \beta;\ \gamma,\ z)$ for fixed $\alpha,\ \gamma,\ z$ and large by module $\beta$ is *erroneously* given as $O((\beta z)^{-\alpha})+O((\beta z)^{\alpha-\gamma}\exp(\beta z))$, $O(\frac{1}{\sqrt{n}}\exp(\frac{n+1}{2}\cdot(1-\frac{1}{4b^2})))$ for our case - if this is correct, we would have a disproof of the Riemann hypothesis! Correct asymptotic for our case can be derived e.g. from the discussion on p. 241 of [13], where the author closely follows 1918 year paper of Watson [14] to obtain, for our case, asymptotic $O(\frac{1}{\sqrt{n}}(2b)^n)$. By these reasons we did not rely on these general results and instead present here Lemma 3.

Now we prove theorem 5. For any fixed $n$, in the limit of large $b$, all terms in the r.h.s. of (7) are $O(1)$ or smaller, either trivially or by the use of Lemma 2, apart from the *positive* $O(\ln b)$ term $\frac{n}{2}\psi\left(\frac{b+1}{2}\right)$ and the derivative $\frac{1}{(n-1)!}\frac{d^n}{dz^n}((z+b)^{n-1}\ln(\varsigma(z)))|_{z=1+b}$ which order is unknown in advance. To evaluate this derivative, we replace it with an integral as in eq. (8), put *c=b>1*, and apply Lemma 3. Then we exploit a well-known fact that for any $b \geq 1+\delta$, where $\delta$ is some arbitrary small fixed positive constant, $|\ln(\varsigma(b+it))|$ is bounded and thus its module is smaller than some appropriate constant: $|\ln(\varsigma(b+it))| < C_1$.

Thus, for *b*>1, we get $\frac{1}{2\pi}|\int_{-\infty}^{\infty}\ln(\varsigma(b+it))\frac{(2b+it)^{n-1}}{(-1+it)^{n+1}}dt| < \frac{C_1}{2\pi}\int_{-\infty}^{\infty}|\frac{(2b+it)^{n-1}}{(-1+it)^{n+1}}|dt < C_1\frac{2b}{4b^2-1}(2b)^{n+1}$. Clearly, for $b \to +\infty$ the constant $C_1$, which is a bound for $|\ln(\varsigma(b+it))|$, is $O(2^{-b})$ whence the value of $\frac{1}{(n-1)!}\frac{d^n}{dz^n}((z+b)^{n-1}\ln(\varsigma(z)))|_{z=1+b} = -\frac{n}{2\pi}\int_{-\infty}^{\infty}\ln(\varsigma(c+it))\frac{(c+b+it)^{n-1}}{(c-b-1+it)^{n+1}}dt$ is in this limit exponentially small and can be neglected. Thus the expression (7) is dominated by the positive term $\frac{n}{2}\psi\left(\frac{b+1}{2}\right)$ which finishes the proof.

***Remark 3.*** The same Theorem 5 can be obtained starting from the expression for the sums at question given in [1] as "arithmetic interpretation" of the generalized Li's sums:

$$\sum_{\rho}(1-\left(\frac{\rho-a}{\rho+a-1}\right)^n) = \sum_{\rho}(1-\left(\frac{\rho+a-1}{\rho-a}\right)^n) = 2-(-1+\frac{1}{a})^n-(-1-\frac{1}{a-1})^n+$$
$$\sum_{j=1}^{n}C_n^j(2a-1)^j\frac{(-1)^j}{(j-1)!}\sum_{m=1}^{\infty}\frac{\Lambda(m)\ln^{j-1}m}{m^a}+\frac{n}{2}(2a-1)(\psi(a/2)-\ln\pi)+ \qquad (14).$$
$$\sum_{j=2}^{n}C_n^j(-1)^j2^{-j}(2a-1)^j\varsigma(j,\ a/2)$$

Indeed, here one easily recognizes our above contributions from zeroes at *z=1* and *z=0*, contributions from derivatives

$\frac{1}{(n-1)!}\frac{d^n}{dz^n}((z+b)^n \ln(\varsigma(z)) = -\sum_{m=1}^{\infty}\frac{\Lambda(m)}{m^z}L_{n-1}^1((z+b)\cdot \ln m)$ (in our eqs. (9, 10) we made one step further), and so on. (For example, after some transformations we can see that the term $\sum_{j=2}^{n} C_n^j (-1)^j 2^{-j}(2a-1)^j \varsigma(j,\ a/2)$ in eq. (14) expresses the contribution from trivial Riemann function zeroes).

This arithmetic interpretation has been obtained applying so called Explicit Formula of Weil, see [2, 15, 16], which relies onto certain Mellin transforms, and from here we see deep analogies between this same Explicit Formula of Weil and approach of the current paper. Using the same Mellin transform-based approach and developing a reasoning of [1, 2], it is easy to show that if there exists a real number $1/2 < \sigma_0 < 1$, such that for all Riemann function zeroes one has $\mathrm{Re}\,\rho \le \sigma_0 + \varepsilon$ (where $\varepsilon$ is an arbitrary small fixed positive number), the following limit holds: $\frac{\varsigma'(a)}{\varsigma(a)} = -\lim_{N\to\infty}(\sum_{m\le N}\frac{\Lambda(m)}{m^a} - \frac{N^{1-a}}{1-a})$ [17]. Differentiation of this equality with respect to *a* readily gives a number of equalities involving higher order derivatives $\frac{d^j}{ds^j}\frac{\varsigma'(s)}{\varsigma(s)}|_{s=a}$ and sums $\sum_{m\le N}\frac{\Lambda(m)\ln^j m}{m^a}$, which can be substituted into above relations expressing the sums over Riemann function zeroes via certain derivatives involving the logarithm of the Riemann zeta-function. This, together with the eq. (3), readily gives a "conditional" arithmetic interpretation for the case $\mathrm{Re}\, z \ge \sigma_0 + \varepsilon$.

***Remark 4.*** Similar theorem can be established for a number of other Riemann zeta-function. We will not pursue this line of researches here, and would like only to sketch the proof for the Dedekind zeta-function, which for an algebraic number field *k* with $r_1$ real and $r_2$ imaginary places is defined as $\varsigma_k(s) = \prod_p (1 - Np^{-s})^{-1}$ for Re*s*>1; product is over all finite prime divisors of *k*. Li's criterion for this function [3] is straightforwardly generalized to give appropriate generalized criterion, and then we are willing to use the fact that the function $Z_k(s) = G_1(s)^{r_1} G_2(s)^{r_2} \varsigma_k(s)$ with $G_1(s) = \pi^{-s/2}\Gamma(s/2)$ and $G_2(s) = (2\pi)^{1-s}\Gamma(s)$ is analytic in the complex plane except for the simple poles at *s=0* and *s=1*. Hence an entire function $\xi_K(s) = As(s-1)|D|^{s/2} Z_k(s)$, where *D* is the discriminant of *k* and *A* is an irrelevant for our purposes constant can be introduced, and the same reasoning as applied before can be repeated. Quite similarly to the case of $\xi(z) = \frac{1}{2} z(z-1)\pi^{-z/2}\Gamma(z/2)\varsigma(z)$, where the behavior of the generalized sums at large *b* and fixed *n* was dominated by the digamma-function $\psi(b/2)$, such a behavior pertinent to $\xi_K(s) = As(s-1)|D|^{s/2} Z_k(s)$ is dominated by the positive sum of digamma functions of $r_1\psi(b/2)$ and $r_2\psi(b)$ leading to the same conclusions.

## 4. Concluding remarks

Certainly, what is proven in Theorem 5 is far from the proof of the Riemann hypothesis, where the non-negativity in another limit, viz. that of some fixed *b>-1/2* (which can be arbitrary large) and $n \to +\infty$, is required. Still, in our opinion, the context of this theorem is interesting enough and does bring some indirect support to the Riemann hypothesis.

What can be said about this another limit? Analyzing for it all terms of (7), we see that all of them are of the order of *O(n*ln*n)* or smaller (again either trivially or by the use of Lemma 2) apart from the exponentially large term $|\frac{1}{2b+1}(1+(-1)^{n+1}(1+\frac{1}{b})^{n})|$ and the derivative $\frac{1}{(n-1)!}\frac{d^n}{dz^n}((z+b)^{n-1}\ln(\varsigma(z)))|_{z=1+b}$ which order is unknown in advance.

From this, certain asymptotic of this latter (or equal to it integral from eq. (8)) can be easily established and the corresponding result can be formulated as necessary and/or sufficient conditions for the Riemann hypothesis. We, however, will not do this here for the following simple reason. The term $\frac{1}{2b+1}(1+(-1)^{n+1}(1+\frac{1}{b})^{n})$ originates from the factor *z-1* in eq. (3), i.e. from a formal “zero” *z*=1 lying to the right to the line Rez=1/2; exactly such zeroes, as we know from the proof of the generalized Li’s criterion [1], lead to the appearance of exponentially large negative sums $\sum_{\rho}(1-\left(\frac{\rho+b}{\rho-1-b}\right)^{n})$ for some values of *n* and, correspondingly, also to the appearance of negative derivatives $\frac{1}{(n-1)!}\frac{d^n}{dz^n}((z+b)^{n-1}\ln(\xi(z)))|_{z=1+b}$. But this same zero at *z=1* is quite formal; it is killed by a simple pole of the Riemann zeta-function, and there is no doubt in the corresponding compensation of the exponentially large term.

Indeed, this compensation can be easily demonstrated in the “first order approximation” to the prime distribution. For this, let us estimate the sum in the r.h.s. of eq. (10) as an integral $I=-\int_{1}^{\infty}x^{-(b+1)}L_{n-1}^{1}((2b+1)\ln x)dx$ (factor *ln(x)* coming from van Mandgoldt function is compensated by *1/ln(x)* density of primes). We apply the variable charge $y=(2b+1)\ln x$ and get $I = -\frac{1}{2b+1}\int_0^\infty \exp(-\frac{b}{2b+1}y)L_{n-1}^1(y)dy$. From an example N 7.414.6 of [11], $\int_0^\infty e^{-qx}L_n(x)dx = \frac{(q-1)^n}{q^{n+1}}$, we have using integration by parts

$\frac{(q-1)^n}{q^{n+1}} = -\frac{1}{q}e^{-qx}L_n(x)|_0^\infty - \frac{1}{q}\int_0^\infty e^{-qx}L_{n-1}^1(x)dx = \frac{1}{q} - \frac{1}{q}\int_0^\infty e^{-qx}L_{n-1}^1(x)dx$, thus

$$\int_0^\infty e^{-qx}L_{n-1}^1(x)dx = 1 - \left(1 - \frac{1}{q}\right)^n \qquad (20).$$

For $q = \frac{b}{2b+1}$ we clearly obtain $1 + (-1)^{n+1}\left(1 + \frac{1}{b}\right)^n$ here which is, by comparison with the term $\frac{1}{2b+1}(1 + (-1)^{n+1}\left(1 + \frac{1}{b}\right)^n)$ in eq. (7), the required exact compensation.

To avoid this situation, let us consider the term $\frac{1}{(n-1)!}\frac{d^n}{dz^n}((z+b)^{n-1}(z-1)\ln(\varsigma(z)))|_{z=1+b}$ where the expression to be differentiated does not contain neither a pole nor a zero at *z=1*. We have a slight modification of eq. (7):

$$\frac{1}{(n-1)!}\frac{d^n}{dz^n}((z+b)^{n-1}\ln(\xi(z)))|_{z=1+b} =$$
$$-\frac{1}{2b+1}\sum_{k=1}^{\infty}(1 - \left(\frac{2k-b}{2k+b+1}\right)^n - \frac{n(2b+1)}{2k+b+1}) + \frac{n}{2}\psi\left(\frac{b+1}{2}\right) + \frac{n}{b+1} - \frac{n\ln\pi}{2}$$
$$+\frac{1}{(n-1)!}\frac{d^n}{dz^n}((z+b)^{n-1}(z-1)\ln(\varsigma(z)))|_{z=1+b} \qquad (7a).$$

Using Lemma 2, we may now formulate the following (not very interesting) sufficient criterion for the Riemann hypothesis.

**Proposition 1.** Riemann hypothesis holds true if for some real $\varepsilon > 0$ asymptotically we have $\frac{1}{(n-1)!}\frac{d^n}{dz^n}((z+b)^{n-1}(z-1)\ln(\varsigma(z)))|_{z=1+b} \ge -(1-\varepsilon)\frac{n}{2}\ln n$.

Our next aim is to calculate the term $\frac{1}{(n-1)!}\frac{d^n}{dz^n}((z+b)^{n-1}\ln(\xi(z)))|_{z=1+b}$ which is, of course, the following sum over non-trivial Riemann function zeroes $k_{n,b}=\sum_\rho(1-\left(\frac{\rho+b}{\rho-1-b}\right)^n)=\sum_\rho(1-\left(\frac{\rho-b-1}{\rho+b}\right)^n)$, *assuming RH*, and see what can be then said about $\frac{1}{(n-1)!}\frac{d^n}{dz^n}((z+b)^{n-1}(z-1)\ln(\varsigma(z)))|_{z=1+b}$ in these circumstances. The corresponding calculation is done in Appendix, and incorporating these results together with what is shown in Lemma 2, we see almost perfect compensation of all terms in eq. (7a). We get:

$$\frac{1}{(n-1)!}\frac{d^n}{dz^n}((z+b)^{n-1}(z-1)\ln(\varsigma(z)))|_{z=1+b}=\frac{n}{2}\ln(b+1/2)-A-\frac{n}{2}\psi\left(\frac{b+1}{2}\right)-\frac{n}{b+1}+o(n)$$

(21)

where $\frac{n}{2}\ln(2-\frac{5}{b+3})<A<\frac{n}{2}\ln(2-\frac{1}{b+1})$; this is a term coming from the "imperfect" estimation of the corresponding sum in Lemma 2. For large $b$ we have $\frac{n}{2}\psi\left(\frac{b+1}{2}\right)=\frac{n}{2}(\ln(b+1)-\ln 2-\frac{1}{b+1}+O(1/b^2))$ so if $b$ tends to infinity in such a manner that $n/b=o(n)$, one has

$$\frac{1}{(n-1)!}\frac{d^n}{dz^n}((z+b)^{n-1}(z-1)\ln(\varsigma(z)))|_{z=1+b}=o(n).$$

Of course, this "almost perfect compensation" of all terms is not at all by chance, it should be somehow related with the fact of the very existence of an appropriate analytical continuation of the Riemann zeta-function to the whole complex plane. We suspect, but cannot prove, that an exact compensation indeed occurs and $\frac{1}{(n-1)!}\frac{d^n}{dz^n}((z+b)^{n-1}(z-1)\ln(\varsigma(z)))|_{z=1+b}=o(n)$ for all $b$>0.

***Remark 5.*** Naturally, one can apply the generalized Littlewood theorem directly to the Riemann zeta-function. Using information about its trivial zeroes and the same function $g(z)=\frac{n(2b+1)(z+b)^{n-1}}{(z-b-1)^{n+1}}-\frac{n(2b+1)}{(z-b-1)^2}$, we get instead of eq. (21):

$$\frac{1}{(n-1)!}\frac{d^n}{dz^n}((z+b)^{n-1}(z-1)\ln(\varsigma(z)))|_{z=1+b}=\frac{1}{2b+1}\sum_{\rho}(1-\left(\frac{\rho+b}{\rho-1-b}\right)^n-\frac{n(2b+1)}{\rho-1-b})$$
$$+\frac{1}{2b+1}\sum_{k=1}^{\infty}(1-\left(\frac{2k-b}{2k+b+1}\right)^n-\frac{n(2b+1)}{2k+b+1})+n\frac{\varsigma'}{\varsigma}(b+1)+\frac{n}{b+1} \qquad (22).$$

Here again, with the use of Lemma 3 and Appendix, the “almost perfect compensation” assuming RH is seen.

**Appendix. Calculation of an asymptotic of generalized Li’s sums and corresponding derivatives assuming the Riemann Hypothesis**

Let us now calculate asymptotic of the sums $k_{n,b}=\sum_{\rho}(1-\left(\frac{\rho+b}{\rho-1-b}\right)^n)$ over non-trivial Riemann function zeroes for large $n$ (and thus also an asymptotic of related with them derivatives) *assuming the Riemann hypothesis*.

**Theorem 6.** *Assume RH. Then for large enough n, for any real fixed* $b\neq -1/2$

$$k_{n,b}=\sum_{\rho}(1-\left(\frac{\rho+b}{\rho-1-b}\right)^n)=\sum_{\rho}(1-\left(\frac{\rho-b-1}{\rho+b}\right)^n)=$$
$$\frac{|2b+1|}{2}n\ln n+\frac{|2b+1|}{2}(\gamma-1-\ln\frac{2\pi}{|2b+1|})n+o(n) \qquad (A1).$$

***Proof.*** Proof is a straightforward generalization of the method presented in Coffey paper [10] (see also [18]; similar asymptotic for Li’s sums was also obtained by Voros with another method [19]). Let us first put $b > -1/2$. Using $\rho = 1/2 + iT$, we write for an argument $\vartheta$ of the function $\frac{\rho + b}{\rho - b - 1}$:

$$\tan\vartheta = -\frac{(2b+1)T}{T^2 - 1/4 - b - b^2} = -\frac{(2b+1)T}{T^2 - (2b+1)^2/4} \qquad \text{(A2).}$$

Correspondingly, $\sin\vartheta = -\frac{(2b+1)T}{T^2 + (2b+1)^2/4}$ and $\cos\vartheta = \frac{T^2 - (2b+1)^2/4}{T^2 + (2b+1)^2/4}$; here we used $(T^2 - 1/4 - b - b^2)^2 + (2b+1)^2 T^2 = (T^2 + 1/4 + b + b^2)^2$. Derivative $d\vartheta/dT$ is found from (A2): $\frac{d\vartheta}{dT} = \frac{2b+1}{T^2 + (2b+1)^2/4}$, and now we are in a position to calculate the sum at question on RH: $k_{n,b} = \Sigma_\rho (1 - (\frac{\rho + b}{\rho - b - 1})^n) = 2\sum_\rho (1 - \cos(n\vartheta_\rho))$ so that, expressed as an integral over the number of non-trivial zeroes *dN*, $k_{n,b} = 2\int_0^\infty (1 - \cos(n\vartheta(T)))dN$. Integrating by parts, we obtain

$$k_{n,b} = 2\int_0^\infty (1 - \cos(n\vartheta(T)))dN = -2n\int_0^\infty \sin(n\vartheta)\frac{d\vartheta}{dT}N(T)dT \qquad \text{(A3)}$$

and then use the approximations $N(T) = \frac{T}{2\pi}\ln\frac{T}{2\pi} - \frac{T}{2\pi} + O(\ln T)$ [4], $\vartheta = -\frac{2b+1}{T} + O(1/T^3)$, $\frac{d\vartheta}{dT} = \frac{2b+1}{T^2} + O(1/T^4)$ to get $k_{n,b} = 2n\int_{T_1}^\infty \frac{(2b+1)}{T^2}\sin(\frac{(2b+1)n}{T})N(T)dT + o(n)$ where $T_1$=14, say (the first zero lies at ½+*i14.1347…*[4]). With the variable change $y = \frac{(2b+1)n}{T}$, we have further $k_{n,b} = 2\int_0^{(2b+1)n/T_1} \sin(y)N(\frac{(2b+1)n}{y})dy = -\frac{(2b+1)n}{\pi}\int_0^\infty \frac{\sin y}{y}(\ln\frac{2\pi y}{(2b+1)n} + 1)dy + o(n)$, and, using examples N3.721.1 $\int_0^\infty \frac{\sin y}{y}dy = \frac{\pi}{2}$ and N4.421.1 $\int_0^\infty \ln y\frac{\sin y}{y}dy = -\frac{\pi}{2}\gamma$ from GR book [11], finally obtain $k_{n,b}=\sum_{\rho}(1-\left(\frac{\rho+b}{\rho-1-b}\right)^n)=\sum_{\rho}(1-\left(\frac{\rho-b-1}{\rho+b}\right)^n)=$
$\frac{2b+1}{2}n\ln n+\frac{2b+1}{2}(\gamma-1-\ln\frac{2\pi}{(2b+1)})n+o(n)$.

The case *b<-1/2* is quite similar with changes of signs whenever appropriate, and in this manner we recover the equation (A1).

***Remark 6.*** Following [10, 18], the sum (derivative) at question can be rewritten, using $\sin n\vartheta=\sin\vartheta\cdot U_{n-1}(\cos\vartheta)$, where $U_k$ is the $k$-th Chebyshev polynomial of the second kind [20], in a rather elegant form

$$k_{n,b}=2n\int_0^{\infty}\frac{(2b+1)^2T}{(T^2+(2b+1)^2/4)^2}U_{n-1}(\frac{T^2-(2b+1)^2/4}{T^2+(2b+1)^2/4})N(T)dT \qquad \text{(A4).}$$

We will not use any properties of this polynomial below, but would like to note the next logical step which is the variable change $x=\frac{T^2-(2b+1)^2/4}{T^2+(2b+1)^2/4}=1-\frac{(2b+1)^2/2}{T^2+(2b+1)^2/4}$. Clearly, $dx=\frac{T(2b+1)^2}{(T^2+(2b+1)^2/4)^2}dT$ so that

$$k_{n,b}=2n\int_{-1}^{1}U_{n-1}(x)N(x)dx \qquad \text{(A5).}$$

Using $T=\frac{1}{2}(2b+1)\sqrt{\frac{1+x}{1-x}}$ and limiting ourselves with the $N(T)=\frac{T}{2\pi}\ln\frac{T}{2\pi}-\frac{T}{2\pi}+O(\ln T)$ precision, we may write $N(x)=\frac{2b+1}{4\pi}\sqrt{\frac{1+x}{1-x}}\ln(\frac{2b+1}{4\pi}\sqrt{\frac{1+x}{1-x}})-\frac{2b+1}{4\pi}\sqrt{\frac{1+x}{1-x}}+O(\ln\frac{1+x}{1-x})$ which is to be substituted into (A3). Note, that integrals of the type $\int_{-1}^{1}U_n(x)(1-x)^{\alpha}(1+x)^{\beta}dx$

quite naturally appear in some applications of the Chebyshev polynomials; see e. g. example N 7.347.2 of GR book [11].